\documentclass{amsart}

\usepackage{amsmath,amssymb,amsthm}
\usepackage{mathtools}
\usepackage{enumitem}
\usepackage{xcolor}
\usepackage{quiver}
\usepackage{tikz-cd}
\usepackage{hyperref}
\hypersetup{
  colorlinks=true,
  linkcolor={blue!60!black},
  citecolor={blue!60!black},
  urlcolor={blue!60!black}
}
\usepackage{comment}

\newtheorem{theorem}{Theorem}[section]
\newtheorem{lemma}[theorem]{Lemma}

\newtheorem{corollary}[theorem]{Corollary}

\theoremstyle{definition}
\newtheorem{definition}[theorem]{Definition}
\newtheorem{remark}[theorem]{Remark}

\newcommand{\lam}{\lambda}

\newcommand{\SKB}{\mathsf{SKB}}

\newcommand{\id}{\mathrm{id}}

\newenvironment{invisible}{{\noindent\sc \colorbox{yellow}{Invisible:}\;}\color{gray}}{\medskip}
\excludecomment{invisible}

\title{$\mathsf{SKB}$ is not algebraically coherent}
\author{Andrea Sciandra, Vito Volpe}
\date{}

\address{%
\parbox[b]{0.9\linewidth}{Département de Mathématiques, Université Libre de Bruxelles, Boulevard du Triomphe, B-1050
 Bruxelles, Belgium.}}
\email{andrea.sciandra@ulb.be}
\urladdr{\url{www.andreasciandra.com}}

\address{%
 \parbox[b]{0.9\linewidth}{Dipartimento di Matematica “Federigo Enriques”, Universit\'a degli Studi di Milano, Via Saldini 50, 20133 Milano, Italy.}}
\email{vito.volpe@unimi.it}


\begin{document}

\begin{abstract}
We prove that the category $\mathsf{SKB}$ is not algebraically coherent. As a consequence, it is not locally algebraically cartesian closed. Then the same can be deduced for cocommutative Hopf braces.
\end{abstract}

\maketitle

\tableofcontents

\section{Introduction}
The notion of \textit{algebraically coherent category} first appeared in \cite{CGvdL} as an algebraic version of the notion of \textit{coherent category} \cite{PTJ}. The parallelism stems from the fact that a coherent category $\mathcal{C}$ is a regular category such that every change-of-base functor of the codomain fibration $\mathrm{cod}: \mathrm{Arr}(\mathcal{C}) \longrightarrow \mathcal{C}$ is coherent, that is, preserves finite limits and jointly strongly epimorphic pairs of arrows. Analogously, an algebraically coherent category is a finitely complete category $\mathcal{C}$ for which the coherence is required of every change-of-base functor of the fibration of points $\mathrm{cod}: \mathrm{Pt}(\mathcal{C}) \longrightarrow \mathcal{C}$. 

The property of being algebraically coherent for a category $\mathcal{C}$ has some interesting consequences, particularly when the category $\mathcal{C}$ is semi-abelian. Indeed, in this case the category $\mathcal{C}$ satisfies the so-called \textit{Smith is Huq} condition (SH), is \textit{strongly protomodular} and satisfies the \textit{normality of Higgins commutator of normal subobjects} condition (NH). Knowing that a category is not only semi-abelian, but satisfies these additional conditions is crucial for many results in categorical algebra, in particular to applications in (co)homology theory. 

Semi-abelian categories were introduced by Janelidze, Márki and Tholen \cite{JMT}, building on Bourn's earlier notion of protomodular category \cite{Bourn1991}. A finitely complete category with a zero object is protomodular if and only if it satisfies the Split Short Five Lemma. Semi-abelian categories are defined as exact (in the sense of Barr \cite{Barr}), protomodular categories with finite coproducts and a zero object. These categories are suitable to develop several basic aspects of homological algebra of groups and rings \cite{Bourn2001}, as well as an abstract theory of commutators and of ideals \cite{Bourn2000}. 

There are several examples of semi-abelian categories which are algebraically coherent, such as the categories of groups, rings (not necessarily unitary), Lie algebras and cocommutative Hopf algebras over a field \cite{GSV}.

Skew braces were introduced in \cite{GV} generalizing braces \cite{Rump}, in connection with the study of set-theoretic solutions to the Yang-Baxter equation. A skew brace is a triple $(A,+,\circ)$, where $(A,+)$ and $(A,\circ)$ are groups such that the compatibility condition $a \circ(b+c) = a \circ b -a + a\circ c$ holds for all $a,b,c \in A$. Skew braces form a variety of $\Omega$-groups in the sense of Higgins \cite{Higgins}, hence in particular a semi-abelian category. They share many interesting properties with the category of groups. For instance, the category $\SKB$ of skew braces satisfies the (SH) condition \cite{BFP2023} and is strongly protomodular \cite{Bourn2024}. However, in \cite{AlSt}, it was recently shown that it is not action accessible. 

In these notes, we further prove another difference between $\mathsf{SKB}$ and $\mathsf{Grp}$. More precisely, we show that $\SKB$ (as well as braces) is not algebraically coherent. This is done by showing that the \textit{Three Subobjects Lemma} does not hold for a (skew) brace of dimension 16 and for an infinite family of (skew) braces of dimension $p^4$, with $p\not=2$ prime. The same counterexample shows that also (NH) is not satisfied. Since any locally algebraically cartesian closed category in the sense of \cite{BournGray} is algebraically coherent, we deduce that $\mathsf{SKB}$ (as well as braces) is not locally algebraically cartesian closed. Therefore, also the category of cocommutative Hopf braces is not algebraically coherent, hence it is not locally algebraically cartesian closed.

\section{Preliminaries}

We recall the definitions and results that are central to this work.\medskip

\noindent\textbf{Algebraically coherent semi-abelian categories.} We first recall the following definition:

\begin{definition}
For any category $\mathcal{C}$, the category $\mathsf{Pt}(\mathcal{C})$ of \textbf{points} in $\mathcal{C}$ is the category whose:

\begin{itemize}
    \item objects $(A,B,p,s)$ are diagrams of the kind
    \[
      \begin{tikzcd}
        A \arrow[r, shift left, "p"] & B \arrow[l, shift left, "s"]
      \end{tikzcd}
      \quad \text{with } ps = 1_B;
    \]

    \item morphisms $(f,g) \colon (A,B,p,s) \to (A',B',p',s')$ are diagrams of the kind
    \[
    \begin{tikzcd}[ampersand replacement=\&]
	A \&\& {A^\prime} \\
	\\
	B \&\& {B^\prime}
	\arrow["f", from=1-1, to=1-3]
	\arrow["p", shift left, from=1-1, to=3-1]
	\arrow["{p^\prime}", shift left, from=1-3, to=3-3]
	\arrow["s", shift left, from=3-1, to=1-1]
	\arrow["g"', from=3-1, to=3-3]
	\arrow["{s^\prime}", shift left, from=3-3, to=1-3]
    \end{tikzcd}
    \]
    where both the upward and downward directed squares commute:
    \[
      fs = s'g \qquad \text{and} \qquad p'f = g p.
    \]
\end{itemize}
\end{definition}

For any object $B$ of $\mathcal{C}$, one can also consider the full subcategory $\mathsf{Pt}_{B}(\mathcal C)$ of $\mathsf{Pt}(\mathcal{C})$ whose objects $(A,p,s)$ are \emph{points over $B$}, i.e.\ $B$ is the codomain of the split epi $p$.



\medskip

If the category $\mathcal{C}$ has pullbacks, then for any point $(A,p,s)$ over $B$ and any morphism $f \colon B' \to B$ in $\mathcal{C}$, one can form the pullback of the pair $(f,p)$:
\[
\begin{tikzcd}[ampersand replacement=\&]
	{B^\prime} \&\& {B^\prime \times_B A} \&\& A \\
	\\
	\&\& {B^\prime} \&\& B
	\arrow["{s^\prime}"', dashed, from=1-1, to=1-3]
	\arrow["{sf}", curve={height=-24pt}, from=1-1, to=1-5]
	\arrow["{1_{B'}}"', curve={height=12pt}, from=1-1, to=3-3]
	\arrow["{p_2}", from=1-3, to=1-5]
	\arrow["{p_1}"', from=1-3, to=3-3]
	\arrow["p", shift left, from=1-5, to=3-5]
	\arrow["f"', from=3-3, to=3-5]
	\arrow["s", shift left, from=3-5, to=1-5]
\end{tikzcd}
\]
By the universal property of the pullback, there is a unique morphism $s':B'\to B'\times_{B}A$ in $\mathcal{C}$ such that
\[
  p_1s' = 1_{B'} \qquad \text{and} \qquad p_2s' = sf.
\]
This provides a point $(B' \times_B A,\, p_1,\, s')$ over $B'$ and a morphism of points
\[
  (p_2, f) \colon (B' \times_B A,\, B',\, p_1,\, s') \longrightarrow (A,\, B,\, p,\, s),
\]
which extends to a functor, called the \textbf{change-of-base functor}:
\[
  f^* \colon \mathsf{Pt}_{B}(\mathcal{C}) \to \mathsf{Pt}_{B'}(\mathcal{C}).
\]
\begin{invisible}
\begin{definition}
A split short exact sequence is a diagram of the form 
$$
\begin{tikzcd}[ampersand replacement=\&]
	0 \&\& K \&\& A \&\& B \&\& 0
	\arrow[from=1-1, to=1-3]
	\arrow["k", from=1-3, to=1-5]
	\arrow["p", shift left, from=1-5, to=1-7]
	\arrow["s", shift left, from=1-7, to=1-5]
	\arrow[from=1-7, to=1-9]
\end{tikzcd}
$$
where $p \circ s = 1_B$ and $k =\mathrm{ker}(p)$.
\end{definition}

\begin{definition}
A morphism of split short exact sequences is a diagram of the form
\begin{align} \label{sssm}
\begin{tikzcd}[ampersand replacement=\&]
	0 \&\& K \&\& A \&\& B \&\& 0 \\
	\\
	0 \&\& {K^\prime} \&\& {A^\prime} \&\& {B^\prime} \&\& 0
	\arrow[from=1-1, to=1-3]
	\arrow["k", from=1-3, to=1-5]
	\arrow["u"', from=1-3, to=3-3]
	\arrow["p", shift left, from=1-5, to=1-7]
	\arrow["v", from=1-5, to=3-5]
	\arrow["s", shift left, from=1-7, to=1-5]
	\arrow[from=1-7, to=1-9]
	\arrow["w", from=1-7, to=3-7]
	\arrow[from=3-1, to=3-3]
	\arrow["{k^\prime}", from=3-3, to=3-5]
	\arrow["{p^\prime}", shift left, from=3-5, to=3-7]
	\arrow["{s^\prime}", shift left, from=3-7, to=3-5]
	\arrow[from=3-7, to=3-9]
\end{tikzcd}
\end{align}

where both rows are split short exact sequences, $v \circ k=k^\prime \circ u$, $w \circ p= p^\prime \circ v$, and $s^\prime \circ w = v \circ s$.
\end{definition}

\begin{definition}
Let $\mathcal{C}$ be a finitely complete pointed category. We say that the split short five lemma holds in $\mathcal{C}$, when for each morphism of split short exact sequences \eqref{sssm}, if $u$ and $w$ are isomophisms, then $v$ is an isomorphism as well. 
\end{definition}

\begin{definition}
A finitely complete pointed category $\mathcal{C}$ is protomodular if the split short five lemma holds in $\mathcal{C}$.    
\end{definition}

\begin{definition}
A category $\mathcal{C}$ is semi-abelian if it is pointed, Barr-exact, protomodular and it has binary coproducts.
\end{definition}

\begin{definition}[\cite{CGvdL}]
 A pointed protomodular category $\mathcal{C}$ is \emph{algebraically coherent} if,
 for each cospan of monomorphisms of split extensions over $B$
 \[
 \begin{tikzcd}
   X_{1} \ar[r,"k_1"] \ar[d,"u_1"'] &
   A_1 \ar[r,shift left,"p_1"] \ar[d,"v_1"] &
   B \ar[d,"\id_B"] \ar[l,shift left,"i_1"] \\
   X \ar[r,"k"] &
   A \ar[r,shift left,"p"] &
   B \ar[l,shift left,"i"] \\
   X_{2} \ar[r,"k_2"'] \ar[u,"u_2"] &
   A_2 \ar[r,shift left,"p_2"] \ar[u,"v_2"'] &
   B \ar[u,"\id_B"'] \ar[l,shift left,"i_2"]
 \end{tikzcd}
 \]
 if $v_1$ and $v_2$ are jointly strongly epimorphic in $\mathcal{C}$,
 then so are $u_1$ and $u_2$.
 \end{definition}
\end{invisible}
We are now able to recall the following definition:
\begin{definition}[{\cite[Definition 3.1]{CGvdL}}]
 A category with finite limits is called \textit{algebraically coherent} if for each morphism $f:X\to Y$ in $\mathcal{C}$, the change-of-base functor
 \[
 f^{*}:\mathrm{Pt}_{Y}(\mathcal{C})\to\mathrm{Pt}_{X}(\mathcal{C})
 \]
 is coherent.
 \end{definition}

We recall that a finitely complete category $\mathcal{C}$ is \textit{locally algebraically cartesian closed} when, for any $f:X\to Y$ in $\mathcal{C}$, the change-of-base functor $f^{*}:\mathrm{Pt}_{Y}(\mathcal{C})\to\mathrm{Pt}_{X}(\mathcal{C})$ is a left adjoint \cite{BournGray}. If $\mathcal{C}$ is locally algebraically cartesian closed then it is algebraically coherent by \cite[Theorem 4.5]{CGvdL}.

Any algebraically coherent semi-abelian category is strongly protomodular \cite[Theorem 6.24]{CGvdL} and satisfies the "Smith is Huq" condition (SH): two equivalence relations on a given object always centralise each other (i.e.\ commute in the sense of Smith \cite{Smith,Pedicchio}) as soon as their normalisations commute in the sense of Huq \cite{Huq}. Any algebraically coherent semi-abelian category also satisfies the condition (NH) of normality of Higgins commutators of normal subobjects \cite{CigoliPhd,CGvdL2}, see \cite[Theorem 6.18]{CGvdL}: the Higgins commutator of two normal subobjects of a given object is again a normal subobject, so that it coincides with the Huq commutator. 

Hence, in such a category $\mathcal{C}$, the Higgins commutator and the Huq commutator of two normal subojects $X,Y$ of an object $A$ are the same, and are usually denoted by $[X,Y]$. \medskip

\noindent\textbf{Skew braces.} The other central notion of this work is the following one:
\begin{definition}[{\cite[Definition 1.1]{GV}}]
 A \emph{skew brace} is a triple $(A,+,\circ)$ where $A^{+}:=(A,+)$ and $A^{\circ}:=(A,\circ)$
 are groups sharing the same identity element $1$, satisfying
 \begin{equation}\label{eq:brace}
   a \circ (b + c) = a \circ b- a + a \circ c
   \qquad \text{for all}\, a,b,c \in A.
 \end{equation}
We write $a^{-\circ}$ for the inverse of $a$ in $A^{\circ}$. Morphisms of skew braces are morphisms of groups with respect to both operations. 

We denote the category of skew braces and their morphisms by $\SKB$.
\end{definition}

Given a skew brace $(A,+,\circ)$, there is a group homomorphism
 \[
\lambda:A^{\circ}\to\mathrm{Aut}(A^{+}),\quad a\mapsto\lam_a(x) := -a+ a \circ x
\]
as proven in \cite[Corollary 1.10]{GV}. 
One can also define the $\star$-product $a \star b := \lam_a(b)- b$ so that, since $a\circ b=a+\lambda_a(b)$, one has $a \circ b = a + (a \star b) + b$. We recall that a \emph{normal subobject} of $A$ in $\SKB$ -- equivalently, an
 \emph{ideal} of $(A,+,\circ)$ -- is a subset $X \subseteq A$ which is simultaneously
 a normal subgroup of $A^{+}$, a normal subgroup of $A^{\circ}$, and satisfies
 $\lam_a(X) \subseteq X$ for all $a \in A$. 

We recall that, given a skew brace $(A,+,\circ)$ and two ideals $I,J$ of $A$, the Huq commutator $[I,J]$ is defined by
\[
\langle[x,y]_{+},[x,y]_{\circ}, x\star y\ |\ x\in I, y\in J\rangle_{A},
\]
where $[\cdot,\cdot]_{+}$ and $[\cdot,\cdot]_{\circ}$ denote the additive and multiplicative commutators, and $\langle\cdot \rangle_{A}$ denotes the ideal generation in $A$ \cite{BFP2023}, see also \cite{GLV,Tsang}. We point out that $[I,J]$ can also be described as $\langle[x,y]_{+},x\star y,y\star x\ |\ x\in I, y\in J\rangle_{A}$.

\section{Counterexamples}

We recall the following result:

\begin{theorem}[{\cite[Theorem 7.1]{CGvdL}}]\label{thm:algcoh}
If $K,L$ and $M$ are normal subobjects of an object $A$ in an algebraically coherent semi-abelian category, then
\[
[[L,M],K]\leq[[K,L],M]\vee[[M,K],L].
\]
\end{theorem}

Given $L=M$, one has
\begin{equation}\label{eq:conditionAC}
[[M,M],K]\leq[[K,M],M]\vee[[M,K],M]=
[[M,K],M].
\end{equation}

We recall that the category $\SKB$ is semi-abelian \cite{BFP2023}. We now suppose that $\mathsf{SKB}$ is algebraically coherent. \medskip

\noindent\textbf{Counterexample 1}. Given a skew brace $(A,+,\circ)$, it is known that $\mathrm{ker}(\lambda)$ is an ideal of $A$. 
\begin{invisible}
In fact, clearly $\mathrm{ker}(\lambda)$ is a normal subgroup of $A^{\circ}$. Moreover, given $a,b\in\mathrm{ker}(\lambda)$, we have $a+ b=a\circ\lambda_{a^{-\circ}}(b)=a\circ b$, hence $\mathrm{ker}(\lambda)$ is also a normal subgroup of $A^{+}$. Moreover, given $x\in A$ and $a\in\mathrm{ker}(\lambda)$, we have $x\circ a\circ\overline{x}\in\mathrm{ker}(\lambda)$, hence
\[
\lambda_x(a)=x^{-1}+(x\circ a)=x^{-1}+(x\circ a\circ\overline{x}\circ x)=x^{-1}+(x\circ a\circ\overline{x})+ x\in\mathrm{ker}(\lambda).
\]
\end{invisible}
Therefore, setting $M=A$ and $K=\mathrm{ker}(\lambda)$ in \eqref{eq:conditionAC}, we must have:
\begin{equation}\label{eq:ACSKB}
[[A,A],\mathrm{ker}(\lambda)]\subseteq[[A,\mathrm{ker}(\lambda)],A].
\end{equation}

 Consider the skew brace $(A,+,\circ)$ defined in the following way. We have $(A,+)=(\mathbb{Z}_{4}\times\mathbb{Z}_{4},+)$, where the sum is defined component-wise and $(A,\circ)$ is the generalized \textit{quaternion group} of order 16. More precisely, given $(a,c),(b,d)\in A$, we set $\xi:=2a+c$ and $v_{0}:=(0,1)$, $v_{1}:=(2,3)$, $v_{2}:=(3,3)$ and $v_{3}:=(3,1)$. Then, we define
\[
(a,c)\circ(b,d):=(a+(-1)^{c}b+dv_{\xi}^{1},c+dv_{\xi}^{2})
\]
\begin{invisible}
This is generated by $g_{1}:=(1,0)$ and $g_{2}:=(0,1)$ with relations
\[
g_{1}^{\circ4}=g_{2}^{\circ4}=0,\quad g_{2}^{\circ2}=g_{1}^{\circ2},\quad g_{2}\circ g_{1}\circ g_{2}^{-\circ}=(0,2).
\]
\end{invisible}
This is the skew brace (indeed, brace) SmallSkewbrace(16,78) in $\mathsf{GAP}$ \cite{GAP}. The unique non-trivial ideals of $(A,+,\circ)$ are the following:
\[
I_{1}:=\{(0,0),(2,0)\}\subset I_{2}:=\{(0,0),(2,0),(1,2),(3,2)\}\subset I_{3}:=\{(a,0), (a,2)\ |\ a\in\mathbb{Z}_{4}\}.
\]

Given two subsets $X,Y$ of $A$, we call $X\star Y=\{x\star y\ |\ x\in X,\ y\in Y\}$. We obtain the following result, which have been checked using $\mathsf{GAP}$:
\begin{lemma}
The following hold:
\begin{itemize}
    \item[1)] $\mathrm{ker}(\lambda)=I_{2},\quad A\star A=I_{3}$,
    \item[2)] $A\star I_{1}=0$,\quad $I_{1}\star A=0$,
    \item[3)] $A\star I_{2}=I_{1}, \quad I_{2}\star A=0$,
    \item[4)] $I_{2}\star I_{3}=0,\quad I_{3}\star I_{2}=I_{1}$.
\end{itemize}
\end{lemma}
We recall that $(A,+)$ is abelian. Therefore, we have $[A,A]=I_{3}$ and then $[[A,A],\mathrm{ker}(\lambda)]=[I_{3},I_{2}]=I_{1}$. Moreover, we have $[A,I_{2}]=I_{1}$ and then $[[A,\mathrm{ker}(\lambda)],A]=[I_{1},A]=0$. This is in contraddiction with \eqref{eq:ACSKB}. Therefore, the category $\mathsf{SKB}$ is not algebraically coherent.\medskip

We now provide an infinite family of counterexamples of orders $p^{4}$, with $p\not=2$ prime.\medskip

\noindent\textbf{Counterexample 2}. Consider a prime $p\not=2$, and define the skew brace $(\mathbb{F}_p^4,+,\circ)$, where $+$ is the summation defined component-wise and 
$\circ: \mathbb{F}_p^4 \times \mathbb{F}_p^4 \longrightarrow \mathbb{F}_p^4$ is defined by
\[
(a,b,c,d)\circ(a^\prime,b^\prime,c^\prime,d^\prime):=(a+a^\prime, b + b^\prime, c+c^\prime+aa^\prime, d+d^\prime+(2c-a^2)b^\prime).
\]
One can check that $\circ$ is associative and defines a group structure, for any odd prime $p$. 
\begin{invisible}
We check it is associative.
Given three elements $(a,b,c,d),(a^\prime,b^\prime,c^\prime,d^\prime),(a^{\prime\prime},b^{\prime\prime},c^{\prime\prime},d^{\prime\prime})$ in $\mathbb{F}^4_p$, we have
    \begin{align*}
        [(a,b,c,d)\circ(a^\prime,b^\prime,c^\prime,d^\prime)]\circ (a^{\prime\prime},b^{\prime\prime},c^{\prime\prime},d^{\prime\prime}) &= (a+a^\prime, b + b^\prime, c+c^\prime+aa^\prime,\\
        &d+d^\prime+(2c-a^2)b^\prime) \circ (a^{\prime\prime},b^{\prime\prime},c^{\prime\prime},d^{\prime\prime}) \\
        &=(a+a^\prime+a^{\prime\prime},b+b^\prime+b^{\prime\prime}, c+c^\prime+aa^\prime+c^{\prime\prime}+aa^{\prime\prime}+a^\prime a^{\prime\prime},\\
        &d+d^\prime+2cb^\prime -a^2 b^\prime +d^{\prime\prime}+ (2(c+c^\prime +aa^\prime)-(a+a^\prime)^2)b^{\prime\prime})\\
        &= (a+a^\prime+a^{\prime\prime},b+b^\prime+b^{\prime\prime}, c+c^\prime+c^{\prime\prime}+aa^\prime+ aa^{\prime\prime}+a^\prime a^{\prime\prime},\\
        &d+d^\prime+d^{\prime\prime}+2cb^\prime+2cb^{\prime\prime}+2c^\prime b^{\prime\prime}-a^2b^\prime -a^2b^{\prime\prime}-a^{\prime 2}b^{\prime\prime}).\\
        (a,b,c,d)\circ[(a^\prime,b^\prime,c^\prime,d^\prime)\circ (a^{\prime\prime},b^{\prime\prime},c^{\prime\prime},d^{\prime\prime})] &= (a,b,c,d) \circ (a^\prime+ a^{\prime\prime},b^\prime + b^{\prime\prime}, c^\prime+c^{\prime\prime}+a^\prime a^{\prime\prime},\\
        & d^\prime +d^{\prime\prime}+(2c^\prime -a^{\prime 2})b^{\prime\prime}) \\
        &= (a+ a^\prime+ a^{\prime\prime},b+b^\prime + b^{\prime\prime}, c+ c^\prime+c^{\prime\prime}+a^\prime a^{\prime\prime}+aa^\prime + aa^{\prime \prime},\\
        & d+d^\prime+d^{\prime\prime}+2cb^\prime+2cb^{\prime\prime}+2c^\prime b^{\prime\prime}-a^2b^\prime -a^2b^{\prime\prime}-a^{\prime 2}b^{\prime\prime}).
    \end{align*}
\end{invisible}
The neutral element is given by the 4-tuple (0,0,0,0), 
\begin{invisible}
Indeed,
\begin{align*}
(a,b,c,d)\circ(0,0,0,0)&=(a+0,b+0,c+0+a+ 0,d+0+(2c-a^2)+ 0)\\
        &=(a,b,c,d);\\
        (0,0,0,0)\circ (a,b,c,d)&=(0+a,0+b,0+c+0+ a,0+d+(2+ 0-0^2)+ b)\\
        &=(a,b,c,d).
    \end{align*}
\end{invisible}
and the inverse of an element $(a,b,c,d)$ is given by $(-a,-b,a^2-c,-d+(2c-a^2)b)$. 
\begin{invisible}
Indeed,
    \begin{align*}
        (a,b,c,d) \circ (-a,-b,a^2-c,-d+(2c-a^2)b) &= (0,0,c+a^2-c-a^2, d-d+(2c-a^2)b+(2c-a^2)(-b)) \\
        &=(0,0,0,0).\\
        (-a,-b,a^2-c,-d+(2c-a^2)b) \circ (a,b,c,d)&=\\
        &(0,0,a^2-c+c-a^2,-d+(2c-a^2)b+d+(2(a^2-c)-a^2)b)\\
        &=(0,0,0,(2c-a^2)b-(2c-a^2)b)\\
        &=(0,0,0,0).
    \end{align*}
\end{invisible}
We now verify the compatibility condition \eqref{eq:brace}. We have
\begin{align*}
    &(a,b,c,d) \circ [(a^\prime,b^\prime,c^\prime,d^\prime)+(a^{\prime\prime},b^{\prime\prime},c^{\prime\prime},d^{\prime\prime})] \\&= (a,b,c,d) \circ (a^\prime+a^{\prime\prime},b^\prime + b^{\prime\prime},c^\prime +c^{\prime\prime},d^\prime + d^{\prime\prime})\\
    &=(a+a^\prime+a^{\prime\prime},b+b^\prime + b^{\prime\prime},c+c^\prime +c^{\prime\prime}+aa^\prime + aa^{\prime\prime},d+d^\prime +d^{\prime\prime}+(2c-a^2)(b^\prime+b^{\prime\prime})).
\end{align*}
On the other hand, we have
\begin{align*}
   (a,b,c,d) \circ (a^\prime,b^\prime,c^\prime,d^\prime) &= (a+a^\prime, b + b^\prime, c+c^\prime+aa^\prime, d+d^\prime+(2c-a^2)b^\prime) \\
   (a,b,c,d) \circ (a^{\prime\prime},b^{\prime\prime},c^{\prime\prime},d^{\prime\prime})&= (a+a^{\prime\prime}, b + b^{\prime\prime}, c+c^{\prime\prime}+aa^{\prime\prime}, d+d^{\prime\prime}+(2c-a^2)b^{\prime\prime}),
\end{align*}
hence
\[
(a+a^\prime, b + b^\prime, c+c^\prime+aa^\prime, d+d^\prime+(2c-a^2)b^\prime)-(a,b,c,d)= (a^\prime,b^\prime,c^\prime+aa^\prime,d^\prime+(2c-a^2)b^\prime)
\]
and then
\begin{align*}
&(a^\prime,b^\prime,c^\prime+aa^\prime,d^\prime+(2c-a^2)b^\prime) + (a+a^{\prime\prime}, b + b^{\prime\prime}, c+c^{\prime\prime}+aa^{\prime\prime}, d+d^{\prime\prime}+(2c-a^2)b^{\prime\prime}) =\\
&=(a+a^\prime+a^{\prime\prime},b+b^\prime + b^{\prime\prime},c+c^\prime +c^{\prime\prime}+aa^\prime + aa^{\prime\prime},d+d^\prime +d^{\prime\prime}+(2c-a^2)(b^\prime+b^{\prime\prime})).
\end{align*}
Therefore, $(\mathbb{F}_p^4,+,\circ)$ is a skew brace, for every odd prime $p$. \medskip

The $\lambda$-action for a tuple $\mathbf{a}=(a,b,c,d) \in \mathbb{F}^4_p$ over an element $\mathbf{x}=(x_1,x_2,x_3,x_4)$ is given by 
\begin{align*}
    \lambda_\mathbf{a}(\mathbf{x})&=-\mathbf{a}+ \mathbf{a} \circ \mathbf{x} \\
    &=(-a,-b,-c,-d) + (a+x_1,b+x_2,c+x_3+ax_1, d+x_4+(2c-a^2)x_2) \\
    &=(x_1,x_2,x_3+ax_2,x_4+(2c-a^2)x_2).
\end{align*}

Consider the skew brace morphisms
\begin{align*}
    \pi_{(0,1,0,0)}: (\mathbb{F}^4_p,+,\circ) &\longrightarrow (\mathbb{F}^4_p,+,\circ)\\
    (a,b,c,d) &\longmapsto (0,b,0,0)
\end{align*}
\begin{invisible}
    It is clearly a group morphism with respect to +. Let us verify that it is a group morphism with respect to $\circ$. Consider $\mathbf{u}=(a,b,c,d),\mathbf{v}=(a^\prime,b^\prime,c^\prime,d^\prime) \in \mathbb{F}^4_p$ and compute
    \begin{align*}
    \pi(\mathbf{u} \circ \mathbf{v}) &= \pi ((a+a^\prime, b+b^\prime,c+c^\prime+aa^\prime,d+d^\prime+(2c-a^2)b^\prime)) \\
    &= (0, b+b^\prime,0,0).\\
    \pi(\mathbf{u}) \circ \pi(\mathbf{v}) &= (0,b,0,0) \circ (0,b^\prime,0,0) \\
    &= (0, b+b^\prime,0,0) 
    \end{align*}
    and notice that the kernel is $\mathrm{ker}(\pi) = \{(a,b,c,d) \in \mathbb{F}^4_p \mid b=0\}=I$.\\
\end{invisible}
and
\begin{align*}
    \pi_{(1,0,1,0)}: (\mathbb{F}^4_p,+,\circ) &\longrightarrow (\mathbb{F}^4_p,+,\circ)\\
    (a,b,c,d) &\longmapsto (a,0,c,0)
\end{align*}
\begin{invisible}
    It is clearly a group morphism with respect to +. Let us verify that it is a group morphism with respect to $\circ$. Consider $\mathbf{u}=(a,b,c,d),\mathbf{v}=(a^\prime,b^\prime,c^\prime,d^\prime) \in \mathbb{F}^4_p$ and compute
    \begin{align*}
    \pi(\mathbf{u} \circ \mathbf{v}) &= \pi ((a+a^\prime, b+b^\prime,c+c^\prime+aa^\prime,d+d^\prime+(2c-a^2)b^\prime)) \\
    &= (a+a^\prime, 0,c+c^\prime + aa^\prime,0).\\
    \pi(\mathbf{u}) \circ \pi(\mathbf{v}) &= (a,0,c,0) \circ (a^\prime,0,c^\prime,0) \\
    &= (a+a^\prime, 0,c+c^\prime + aa^\prime,0)
    \end{align*}
    and notice that the kernel is $\mathrm{ker}(\pi) = \{(a,b,c,d) \in \mathbb{F}^4_p \mid a=0 \land c=0\}=J$.\\
\end{invisible}
and the ideals 
\[
I = \{(a,0,c,d) \mid a,c,d \in \mathbb{F}_p^4\}=\mathrm{ker}(\pi_{(0,1,0,0)}),\quad J = \{(0,b,0,d) \mid b,d \in \mathbb{F}_p^4\}=\mathrm{ker}(\pi_{(1,0,1,0)}).
\]
We now compute $[I,J]$ on generators. Given $i:=(a,0,c,d) \in I$ and $j:=(0,b,0,d^\prime) \in J$, we have:
\begin{itemize}
    \item $[i,j]_+=(0,0,0,0)$ since $(\mathbb{F}^4_p,+)$ is an abelian group.
    \item $[i,j]_\circ$ is given by:
    \begin{align*}
    [i,j]_\circ &=i\circ j  \circ i^{-\circ} \circ j^{-\circ}   \\
    &=(a,b,c,d+d^\prime+(2c-a^2)b)\circ (-a,0,a^2-c,-d)\circ j^{-\circ} \\
    &=(0,b,0,d^\prime +(2c-a^2)b) \circ (0,-b,0,-d^\prime)\\
    &=(0,0,0,(2c-a^2)b).
    \end{align*}
    \item $i*j:=-i+i \circ j - j$ is given by:
    \begin{align*}
    i*j &= (-a,0,-c,-d)+(a,b,c,d+d^\prime+(2c-a^2)b)+(0,-b,0,-d^\prime) \\
    &= (0,0,0,(2c-a^2)b) = [i,j]_\circ .
    \end{align*}
\end{itemize}
One can easily see that $\{(0,0,0,a) \mid a \in \mathbb{F}^4_p\}$ is an ideal, hence it coincides with $[I,J]$. This is clearly central with respect to both the two operations and $x+y=x\circ y$ with $x\in [I,J]$ and $y\in\mathbb{F}^4_p$. Therefore, we get $[[I,J],I]=\{0\}$.

Now, we compute $[I,I]$. Given $i=(a,0,c,d), j=(a^\prime,0,c^\prime,d^\prime) \in I$, we have:
\begin{itemize}
    \item $[i,j]_+=(0,0,0,0)$ since $(\mathbb{F}^4_p,+)$ is an abelian group.
    \item $[i,j]_\circ$ is given by:
    \begin{align*}
        [i,j]_\circ &= (a+a^\prime,0,c+c^\prime +aa^\prime, d+d^\prime) \circ (-a,0,a^2-c,-d) \circ j^{-\circ}\\
        &= (a^\prime,0,c^\prime,d^\prime) \circ (-a^\prime,0,a^{\prime 2}-c^\prime,-d^\prime) \\
        &= (0,0,0, 0).
    \end{align*}
    \item $i*j$ is given by:
    \begin{align*}
        i*j&=(-a,0,-c,-d)+(a+a^\prime,0,c+c^\prime +aa^\prime, d+d^\prime)+ (-a^\prime,0,-c^\prime,-d^\prime)\\
        &=(0,0,aa^\prime, 0)
    \end{align*}
\end{itemize}
Considering the ideal generated by $\{(0,0,a,0) \mid a \in \mathbb{F}_p\}$, we get $[I,I]=\{(0,0,x,-2xy)\ |\ x,y\in\mathbb{F}^{4}_{p}\}$. We now compute $[[I,I],J]$. Given $i=(0,0,x,-2xy) \in [I,I]$ and $j=(0,b,0,d) \in J$, we compute:
\begin{itemize}
    \item $[i,j]_+=0$ as above;
    \item $[i,j]_\circ$ is given by:
    \begin{align*}
        [i,j]_\circ &= (0,b,x,d-2xy+2xb) \circ (0,0,-x,2xy) \circ j^{-\circ}\\
        &=(0,b,0,d+2xb) \circ (0,-b,0,-d) \\
        &= (0,0,0,2xb).
    \end{align*}
    \item $i*j$ is given by:
    \begin{align*}
        i*j&=(0,0,-x,2xy)+(0,b,x,d-2xy+2xb)+(0,-b,0,-d)\\
        &=(0,0,0,2xb).
    \end{align*}
\end{itemize}

Therefore, we get 
$$
[[I,I],J]=\langle(0,0,0,1)\rangle\neq \{0\},
$$
which again contradicts \eqref{eq:conditionAC}.

\begin{remark}
We notice that for $p=2$ this counterexample fails, since we have $[[I,I],J]=\{0\}$. The case $p=3$, of dimension 81, could be checked using $\mathsf{GAP}$. 
\end{remark}

\begin{remark}
Since our examples are in particular braces, i.e.\ skew braces with $+$ abelian, the semi-abelian subvariety of braces is not algebraically coherent, hence it is not locally algebraically cartesian closed. 
\end{remark}

We observe that cocommutative Hopf braces, Hopf-theoretic generalizations of skew braces introduced in \cite{AGV}, form a semi-abelian category as proven in \cite{GranSciandra}. Moreover, $\mathsf{SKB}$ is a \textit{Birkhoff subcategory} of the category $\mathsf{HBR}_{\mathrm{coc}}$ of cocommutative Hopf braces, that is, a full
replete reflective subcategory that is closed under subobjects and regular quotients \cite[Corollary 6.11]{GranSciandra}. In fact, it is further a localization of $\mathsf{HBR}_{\mathrm{coc}}$, see \cite[Theorem 6.12]{GranSciandra}. Therefore, by \cite[Proposition 3.7]{CGvdL}, we get the following result:

\begin{corollary}
    The category $\mathsf{HBR}_{\mathrm{coc}}$ is not algebraically coherent. Hence it is not locally algebraically cartesian closed.
\end{corollary}

\section{Failure of (NH)}
Recall that a skew brace $(A,+,\circ)$ is an $\Omega$-group (in the sense of \cite{Higgins}). Indeed, $(A,+)$ is a group and $\Omega =\{\circ,(-)^{-\circ}\}$ contains the binary operation $\circ$ and the unary operation $(-)^{-\circ}$. By \cite[Theorem 4B]{Higgins}, given two ideals $I,J$ of the skew brace $(A,+,\circ)$, the Higgins commutator $[I,J]^H$ is the ideal of $I \lor J$ generated by: 
\begin{itemize}
    \item $-i-j+i+j$, for all $i$ in $I$ and $j$ in $J$;
    \item $-i^{-\circ}-j^{-\circ}+(i+j)^{- \circ}$, for all $i$ in $I$ and $j$ in $J$;
    \item $-i_1 \circ i_2 - j_1 \circ j_2 + (i_1 +j_1)\circ (i_2 +j_2) $, for all $i_1,i_2$ in $I$ and $j_1,j_2$ in $J$.
\end{itemize}

Here we provide an example of an ideal $I$ of a skew brace $(A,+,\circ)$ such that $[I,I]^H$ is not an ideal of $(A,+,\circ)$. \medskip

We consider again the skew brace $(\mathbb{F}_p^4,+,\circ)$ of Counterexample 2. and the ideal
\[
I=\{(a,0,c,d) \mid a,c,d \in \mathbb{F}_p\}.
\]
We compute the Higgins commutator $[I,I]^H$ and prove that this is not an ideal of $(\mathbb{F}_p^4,+,\circ)$.

Consider $i=(a,0,c,d), j=(a^\prime,0,c^\prime,d^\prime) \in I$. We now compute: 
\begin{itemize}
    \item $-i-j+i+j = 0$ since $(\mathbb{F}_p^4,+)$ is an abelian group.
    \item $ -i^{-\circ}-j^{-\circ}+(i+j)^{- \circ}$ is given by:
    \begin{align*}
     &-(-a,0,a^2-c,-d) - (-a^\prime,0,a^{\prime 2}-c^\prime, -d^\prime) + (a+a^\prime, 0, c+c^\prime,d+d^\prime)^{-\circ}\\
    &= (a+a^\prime,0,c+c^\prime-a^2-a^{\prime 2},d+d^\prime) + (-a-a^\prime,0,a^2+2aa^\prime+a^{\prime 2} -c-c^\prime,-d-d^\prime)\\
    &=(0,0,2aa^\prime,0).
    \end{align*}
\end{itemize}
Now, given $i_1=(a,0,c,d), i_2=(a^\prime,0,c^\prime,d^\prime), j_1=(e,0,g,h),j_2=(e^\prime,0,g^\prime,h^\prime) \in I$, we compute $-(i_1 \circ i_2) - (j_1 \circ j_2) + (i_1 +j_1)\circ (i_2 +j_2)$. We do it piece by piece 
$$
-(i_1\circ i_2) = -(a+a^\prime,0,c+c^\prime+aa^\prime,d+d^\prime).
$$
$$
-(j_1 \circ j_2) = -(e+e^\prime, 0, g+g^\prime + ee^\prime, h+h^\prime).
$$
Lastly,
\begin{align*}
(i_1 +j_1)\circ (i_2 +j_2) &= (a+e,0,c+g,d+h) \circ (a^\prime + e^\prime,0,c^\prime +g^\prime,d^\prime+h^\prime)\\
&=(a+a^\prime+e+e^\prime,0,c+c^\prime+g+g^\prime +aa^\prime +ae^\prime+a^\prime e+e e^\prime, d+d^\prime+h+h^\prime).
\end{align*}
Putting all together, we get 
\begin{align*}
-(i_1 \circ i_2) - (j_1 \circ j_2) + (i_1 +j_1)\circ (i_2 +j_2) = (0,0,ae^\prime+a^\prime e, 0).
\end{align*}
Since $\{ (0,0,c,0) \mid c \in \mathbb{F}_p^4\}$ is normal in $I = I \lor I$, it coincides with $[I,I]^H$. As was already observed before, this is not an ideal of $(\mathbb{F}_p^4,+,\circ)$, since it is not normal with respect to $\circ$. Take for example $p=3$, $h=(0,0,1,0)$ and $g=(0,1,0,0)$. We have that $g^{-\circ}=(0,-1,0,0)=(0,2,0,0)$. We compute
\begin{align*}
    g\circ h\circ g^{-\circ} &=(0,1,0,0) \circ (0,0,1,0) \circ (0,2,0,0)\\
    &=(0,1,1,0) \circ (0,2,0,0)=(0,3,1, 0+0+(2-0)\cdot 2) = (0,0,1,1) \not\in [I,I]^H.
\end{align*}

\begin{remark}
   As a consequence of \cite[Proposition 2.5]{GvdL}, we further notice that $\SKB$ is not a peri-abelian category \cite[Definition 4.1]{Bourn2010}.
\end{remark}

\end{document}